\documentclass[reqno]{amsart}
\usepackage{amsthm, amssymb, amsfonts}
\usepackage{lmodern}
\usepackage{color} 
\usepackage{hyperref}
\usepackage[T5]{fontenc}
\usepackage{amscd,amssymb}
\usepackage[v2,cmtip]{xy}
\usepackage{mathrsfs}
\theoremstyle{plain}
\newtheorem{thm}{Theorem}[section]
\newtheorem{prop}[thm]{Proposition}

\newtheorem{con}[thm]{Conjecture}
\newtheorem{corl}[thm]{Corollary}
\theoremstyle{definition}
\newtheorem{defn}[thm]{Definition}
\newtheorem{rem}[thm]{Remark}
\newtheorem{nota}[thm]{Notation}

\theoremstyle{plain}

\theoremstyle{definition}

\begin{document} 

\title[A note on the hit problem for the Steenrod algebra]{A note on the hit problem for the Steenrod algebra and its applications} 

\author{Nguyen Khac Tin} 

\address{Faculty of Applied Sciences, 
Ho Chi Minh City University of Technology and Education,
01 V\~o V\u an Ng\^an, Th\h u \DJ\'\uhorn c, H\`\ocircumflex\ Ch\'i Minh city, Viet Nam.}

\email{tinnk@hcmute.edu.vn}

\maketitle

\begin{abstract}
Let $P_{k}=H^{*}((\mathbb{R}P^{\infty})^{k})$ be the modulo-$2$ cohomology algebra of the direct product of $k$ copies of infinite dimensional real projective spaces $\mathbb{R}P^{\infty}$. Then, $P_{k}$ is isomorphic to the graded polynomial algebra $\mathbb{F}_{2}[x_{1},\ldots,x_{k}]$ of $k$ variables, in which each $x_{j}$ is of degree 1, and let $GL_k$ be the general linear group over the prime field $\mathbb{F}_2$ which acts naturally on $P_k$. Here the cohomology is taken with coefficients in the prime field $\mathbb F_2$ of two elements. We study the {\it hit problem}, set up by Frank Peterson, of finding a minimal set of generators for the polynomial algebra $P_k$ as a module over the  mod-2 Steenrod algebra, $\mathcal{A}$. 

In this Note, we explicitly compute the hit problem for $k = 5$ and the degree $5(2^s-1)+24.2^s$ with $s$ an arbitrary non-negative integer. These results are used to study the Singer algebraic transfer which is a homomorphism from the homology of  the mod-$2$ Steenrod algebra, $\mbox{Tor}^{\mathcal{A}}_{k, k+n}(\mathbb{F}_2, \mathbb{F}_2),$ to the subspace of  $\mathbb{F}_2\otimes_{\mathcal{A}}P_k$ consisting of all the $GL_k$-invariant classes of degree $n.$ We show that Singer's conjecture for the algebraic transfer is true in the case $k=5$ and the above degrees. This method is different from that of Singer in studying the image of the algebraic transfer.
Moreover, as a consequence, we get the dimension results for  polynomial algebra in some generic degrees in the case $k=6.$
\end{abstract}

\medskip
\section{\bf Introduction}\label{s1} 
\setcounter{equation}{0}

Denote by $P_{k}=H^{*}((\mathbb{R}P^{\infty})^{k})$ the modulo-$2$ cohomology algebra of the direct product of $k$ copies of infinite dimensional real projective spaces $\mathbb{R}P^{\infty}$. Then, $P_{k}$ is isomorphic to the graded polynomial algebra $\mathbb{F}_{2}[x_{1},x_{2},\ldots,x_{k}]$ of $k$ variables, in which each $x_{j}$ is of degree 1. Here the cohomology is taken with coefficients in the prime field $\mathbb F_2$ of two elements. 

Being the cohomology of a group,  $P_k$ is a module over the mod-2 Steenrod algebra, $\mathcal A$.  The action of $\mathcal A$ on $P_k$ is determined by the elementary properties of the Steenrod squares $Sq^i$ and the Cartan formula (see Steenrod and Epstein~\cite{s-e62}).

A polynomial $f$ in $P_k$  is called {\it hit} if it can be written as a finite sum $f = \sum_{u\geqslant 0}Sq^{2^u}(h_u)$ 
for suitable polynomials $h_u$.  That means $f$ belongs to  $\mathcal{A}^+P_k$,  where $\mathcal{A}^+$ denotes the augmentation ideal in $\mathcal{A}$.

Let $GL_k$ be the general linear group over the field $\mathbb F_2$. This group acts naturally on $P_k$ by matrix substitution. Since the two actions of $\mathcal A$ and $GL_k$ upon $P_k$ commute with each other, there is an action of $GL_k$ on $\mathbb F_2\otimes_{\mathcal A}P_k.$ 

Many authors  study the {\it hit problem} of determination of a minimal set  of generators for $P_k$ as a module over the Steenrod algebra, or equivalently, a basis of $\mathbb{F}_{2}{\otimes}_{\mathcal{A}}P_{k}$. This problem has first been studied by Peterson \cite{pe87}, Wood \cite{wo89}, Singer \cite {si89}, Priddy \cite{pr90},  who show its relationship to several classical problems in homotopy theory. 
Then, this problem was investigated by Nam~\cite{na04}, Silverman~\cite{sil95}, Wood~\cite{wo89}, Sum~\cite{su10, su15, su19}, Tin-Sum~\cite{t-s16}, Tin~\cite{ti2011, ti2026, ti21AEJM} and others.

For a positive integer $n$, by $\mu(n)$ one means the smallest number $r$ for which it is possible to write $n = \sum_{1\leqslant i\leqslant r}(2^{u_i}-1),$ where $u_i >0$. This result implies a result of  Wood, which originally is a conjecture of Peterson~\cite{pe87}.
 
\begin{thm}[Wood~\cite{wo89}]\label{dlwo}
If $\mu(n) > k$, then $(\mathbb{F}_{2}{\otimes}_{\mathcal{A}}P_{k})_n = 0$.
\end{thm} 

From the above result of Wood, the hit problem is reduced to the case of degree $n$ with $\mu(n) \leqslant k.$

Let $GL_k$ be the general linear group over the field $\mathbb F_2$. This group acts naturally on $P_k$ by matrix substitution. Since the two actions of $\mathcal A$ and $GL_k$ upon $P_k$ commute with each other, there is an action of $GL_k$ on $\mathbb{F}_{2}{\otimes}_{\mathcal{A}}P_{k}$. 

For a nonnegative integer $d$, denote by $(P_k)_d$ the subspace of $P_k$ consisting of all the homogeneous polynomials of degree $d$ in $P_k$ and by $(\mathbb{F}_{2}{\otimes}_{\mathcal{A}}P_{k})_d$ the subspace of $QP_k$ consisting of all the classes represented by the elements in $(P_k)_d$. 
One of the major applications of hit problem is in surveying a homomorphism introduced by Singer. In \cite{si89}, Singer defined the algebraic transfer,  which is a homomorphism
$$\varphi_k :\text{Tor}^{\mathcal A}_{k,k+d} (\mathbb F_2,\mathbb F_2) \longrightarrow  (\mathbb{F}_{2}{\otimes}_{\mathcal{A}}P_{k})_d^{GL_k}$$
from the homology of the Steenrod algebra to the subspace of $(\mathbb{F}_{2}{\otimes}_{\mathcal{A}}P_{k})_d$ consisting of all the $GL_k$-invariant classes. It is a useful tool in describing the homology groups of the Steenrod algebra, $\text{Tor}^{\mathcal A}_{k,k+d} (\mathbb F_2,\mathbb F_2)$. The hit problem and the algebraic transfer were studied by many authors (see Boardman ~\cite{bo93}, H\uhorn ng ~\cite{hu05}, Ch\ohorn n and H\`a~ \cite{c-h10}, Nam~ \cite{na04}, Phuc \cite{ph20}, Phuc and Sum \cite{p-s17}, Sum~\cite{su10, su13, su15}, Sum and T\'in~\cite{s-t21} and others).

Singer showed in \cite{si89} that $\varphi_k$ is an isomorphism for $k=1,2$. Boardman showed in \cite{bo93} that $\varphi_3$ is also an isomorphism. However, for any $k\geqslant 4$,  $\varphi_k$ is not a monomorphism in infinitely many degrees (see Singer \cite{si89}, H\uhorn ng \cite{hu05}.) Singer made the following conjecture.

\begin{con}[Singer \cite{si89}] The algebraic transfer $\varphi_k$ is an epimorphism for any $k \geqslant 0$.
\end{con}

 The conjecture is true for $k\leqslant 3$. However, for $k > 3$, it is still open. Recently, the hit problem and its applications to representations of general linear groups have been presented in the books of Walker and Wood~\cite{w-w18, w-w182}.

One of the extremely useful tools for computing the hit problem and studying Singer's transfer is the Kameko squaring operation

$$\widetilde{Sq}^0_*:= (\widetilde{Sq}^0_*)_{(k, d)}: (\mathbb{F}_{2}{\otimes}_{\mathcal{A}}P_{k})_{2k+d} \to (\mathbb{F}_{2}{\otimes}_{\mathcal{A}}P_{k})_{d},$$
which is induced by an $\mathbb F_2$-linear map $\varphi_k: P_k \to P_k$, given by
$$
\varphi_k(x) = 
\begin{cases}y, &\text{if }x=x_1x_2\ldots x_ky^2,\\  
0, & \text{otherwise,} \end{cases}
$$
for any monomial $x \in P_k$. The map $\varphi_k$ is not an $\mathcal A$-homomorphism. However, 
$\varphi_k Sq^{2i} = Sq^{i}\varphi_k$ and $\varphi_k Sq^{2i+1} = 0$ for any non-negative integer $i$.

\begin{thm}[Kameko~\cite{ka90}]\label{dlk} 
Let $d$ be a non-negative integer. If $\mu(2d+k)=k$, then 
$$ (\widetilde{Sq}^0_*)_{(k, d)}: (\mathbb{F}_{2}{\otimes}_{\mathcal{A}}P_{k})_{2d+k}\longrightarrow (\mathbb{F}_{2}{\otimes}_{\mathcal{A}}P_{k})_d$$
is an isomorphism of $GL_k$-modules. 
\end{thm}
Thus, the hit problem is reduced to the case of degree $n$ of the form $$ n=a(2^s-1)+2^sb, $$ where $a, b, m$ are non-negative intergers such that $0 \leqslant \mu(b)<a \leqslant k.$

\medskip
Now, the $\mathbb F_2$-vector space $\mathbb{F}_{2}{\otimes}_{\mathcal{A}}P_{k}$ was explicitly calculated by Peterson~\cite{pe87} for $k=1, 2,$ by Kameko~\cite{ka90} for $k=3$ and by Sum~\cite{su15} for $k = 4$. However, for $k > 4$, it is still unsolved, even in the case of $k=5$ with the help of computers.

For $a=k-1=4$ and $b=0,$ the vector space $(\mathbb{F}_{2}{\otimes}_{\mathcal{A}}P_{5})_n$ is explicitly computed by Phuc and Sum~\cite{p-s17}, and in the case $a=k-2=3, b=1$ by Phuc~\cite{ph18} and by the present author \cite{ti21PJA} in the case $a=k-1=4, \ b \in \{ 8; 10; 11 \}$ with $s$ an arbitrary non-negative integer.

In this Note, we explicitly compute the hit problem for $k = 5$ and the degree $5(2^s-1)+24.2^s$ with $s$ an arbitrary non-negative integer. Using this result, we show that Singer's conjecture for the algebraic transfer is true in the case $k=5$ and the above degrees. Moreover, as a consequence, we get the dimension results for  polynomial algebra in some generic degrees in the case $k=6.$


\section{\bf Preliminaries}\label{s2}
\setcounter{equation}{0}

In this section, we recall some needed information from Kameko~\cite{ka90}, Singer~\cite{si89} and Sum~\cite{su10}, which will be used in the next section.
\begin{nota} We denote $\mathbb N_k = \{1,2, \ldots , k\}$ and
\begin{align*}
X_{\mathbb J} = X_{\{j_1,j_2,\ldots , j_s\}} =
 \prod_{j\in \mathbb N_k\setminus \mathbb J}x_j , \ \ \mathbb J = \{j_1,j_2,\ldots , j_s\}\subset \mathbb N_k,
\end{align*}

Let $\alpha_i(a)$ denote the $i$-th coefficient  in dyadic expansion of a non-negative integer $a$. That means
$a= \alpha_0(a)2^0+\alpha_1(a)2^1+\alpha_2(a)2^2+ \ldots ,$ for $ \alpha_i(a) =0$ or 1 with $i\geqslant 0$. 

Let $x=x_1^{a_1}x_2^{a_2}\ldots x_k^{a_k} \in P_k$. Denote $\nu_j(x) = a_j, 1 \leqslant j \leqslant k$.  
Set 
$$\mathbb J_t(x) = \{j \in \mathbb N_k :\alpha_t(\nu_j(x)) =0\},$$
for $t\geqslant 0$. Then, we have
$x = \prod_{t\geqslant 0}X_{\mathbb J_t(x)}^{2^t}.$ 
\end{nota}
\begin{defn}
For a monomial  $x$ in $P_k$,  define two sequences associated with $x$ by
\begin{align*} 
\omega(x)=(\omega_1(x),\omega_2(x),\ldots , \omega_i(x), \ldots),\ \
\sigma(x) = (\nu_1(x),\nu_2(x),\ldots ,\nu_k(x)),
\end{align*} 
where
$\omega_i(x) = \sum_{1\leqslant j \leqslant k} \alpha_{i-1}(\nu_j(x))= \deg X_{\mathbb J_{i-1}(x)},\ i \geqslant 1.$
The sequence $\omega(x)$ is called  the weight vector of $x$. 

Let $\omega=(\omega_1,\omega_2,\ldots , \omega_i, \ldots)$ be a sequence of non-negative integers.  The sequence $\omega$ is called  the weight vector if $\omega_i = 0$ for $i \gg 0$.
\end{defn}

The sets of all the weight vectors and the exponent vectors are given the left lexicographical order. 

  For a  weight vector $\omega$,  we define $\deg \omega = \sum_{i > 0}2^{i-1}\omega_i$.  
Denote by   $P_k(\omega)$ the subspace of $P_k$ spanned by all monomials $y$ such that
$\deg y = \deg \omega$, $\omega(y) \leqslant \omega$, and by $P_k^-(\omega)$ the subspace of $P_k$ spanned by all monomials $y \in P_k(\omega)$  such that $\omega(y) < \omega$. 

\begin{defn}\label{dfn2} Let $\omega$ be a weight vector and $f, g$ two polynomials  of the same degree in $P_k$. 

i) $f \equiv g$ if and only if $f - g \in \mathcal A^+P_k$. If $f \equiv 0$ then $f$ is called {\it hit}.

ii) $f \equiv_{\omega} g$ if and only if $f - g \in \mathcal A^+P_k+P_k^-(\omega)$. 
\end{defn}

Obviously, the relations $\equiv$ and $\equiv_{\omega}$ are equivalence ones. Denote by $QP_k(\omega)$ the quotient of $P_k(\omega)$ by the equivalence relation $\equiv_\omega$. Then, we have 
$$QP_k(\omega)= P_k(\omega)/ ((\mathcal A^+P_k\cap P_k(\omega))+P_k^-(\omega)).$$  
For a  polynomial $f \in  P_k$, we denote by $[f]$ the class in $QP_k$ represented by $f$. If  $\omega$ is a weight vector, then denote by $[f]_\omega$ the class represented by $f$. Denote by $|S|$ the cardinal of a set $S$.

\begin{defn}\label{defn3} 
Let $x, y$ be monomials of the same degree in $P_k$. We say that $x <y$ if and only if one of the following holds:  

i) $\omega (x) < \omega(y)$;

ii) $\omega (x) = \omega(y)$ and $\sigma(x) < \sigma(y).$
\end{defn}

\begin{defn}
A monomial $x$ is said to be inadmissible if there exist monomials $y_1,y_2,\ldots, y_m$ such that $y_t<x$ for $t=1,2,\ldots , m$ and $x - \sum_{t=1}^my_t \in \mathcal A^+P_k.$ 

A monomial $x$ is said to be admissible if it is not inadmissible.
\end{defn}

Obviously, the set of all the admissible monomials of degree $n$ in $P_k$ is a minimal set of $\mathcal{A}$-generators for $P_k$ in degree $n$. 

Now, we recall some notations and definitions in \cite{su15}, which will be used in the next sections. We set 
\begin{align*} 
P_k^0 &=\langle\{x=x_1^{a_1}x_2^{a_2}\ldots x_k^{a_k} \ : \ a_1a_2\ldots a_k=0\}\rangle,
\\ P_k^+ &= \langle\{x=x_1^{a_1}x_2^{a_2}\ldots x_k^{a_k} \ : \ a_1a_2\ldots a_k>0\}\rangle. 
\end{align*}

It is easy to see that $P_k^0$ and $P_k^+$ are the $\mathcal{A}$-submodules of $P_k$. Furthermore, we have the following.

\begin{prop}\label{2.7} We have a direct summand decomposition of the $\mathbb F_2$-vector spaces
$\mathbb{F}_{2}{\otimes}_{\mathcal{A}}P_{k} =QP_k^0 \oplus  QP_k^+.$
Here $QP_k^0 = \mathbb F_2\otimes_{\mathcal A}P_k^0$ and  $QP_k^+ = \mathbb F_2\otimes_{\mathcal A}P_k^+$.
\end{prop}

\begin{defn} For any $1 \leqslant i \leqslant k$, define the homomorphism $f_i: P_{k-1} \to P_k$ of algebras by substituting
$$f_i(x_j) = \begin{cases} x_j, &\text{ if } 1 \leqslant j <i,\\
x_{j+1}, &\text{ if } i \leqslant j <k.
\end{cases}$$
Then, $f_i$ is a homomorphism of $\mathcal A$-modules.
\end{defn}

For a subset $B \subset P_k,$ we denote $[B] = \{[f]:\,f\in  B \}.$ Obviously, we have
\begin{prop}\label{mdp0}
 It is easy to see that if $B$ is a minimal set of generators for $\mathcal A$-module $P_{k-1}$ in degree $n$, then  $f(B)= \bigcup _{i=1}^kf_i(B)$  is a minimal set of generators for $\mathcal A$-module $P_k^0$ in degree $n$.
\end{prop}

From now on, we denote by $B_{k}(n)$ the set of all admissible monomials of degree $n$  in $P_k$, $B_{k}^0(n) = B_{k}(n)\cap P_k^0$, $B_{k}^+(n) = B_{k}(n)\cap P_k^+$. For a weight vector $\omega$ of degree $n$, we set $B_k(\omega) = B_{k}(n)\cap P_k(\omega)$, $B_k^+(\omega) = B_{k}^+(n)\cap P_k(\omega)$. 

Then, $[B_k(\omega)]_\omega$ and $[B_k^+(\omega)]_\omega$, are respectively the basses of the $\mathbb F_2$-vector spaces $QP_k(\omega)$ and $QP_k^+(\omega) := QP_k(\omega)\cap QP_k^+$.

\section{\bf The Main Results}\label{s3}
\setcounter{equation}{0}

\medskip
We first recall a result in \cite{t-s16} the following: Let $d$ be an arbitrary non-negative integer. Set\ $t(k,d) = \max\{0,k- \alpha(d+k) -\zeta(d+k)\},$ where $\zeta(n)$ the greatest integer $u$ such that $n$ is divisible by $2^u,$ that means $n = 2^{\zeta(n)}m,$ with $m$ an odd integer. 

\medskip
\begin{thm}[Tin-Sum~\cite{t-s16}]\label{dlt-s16} Let $d$ be an arbitrary non-negative integer. Then
$$(\widetilde{Sq}^0_*)^{s-t}: (\mathbb F_2 \otimes_{\mathcal A} P_k)_{k(2^s-1) + 2^sd} \longrightarrow (\mathbb F_2 \otimes_{\mathcal A} P_k)_{k(2^t-1) + 2^td}$$
is an isomorphism of $GL_k$-modules for every $s \geqslant t$ if and only if $t \geqslant  t(k,d).$
\end{thm}

\medskip
It is easy to check that for $k=5$ and $d=53$ then $$t(k,d) = \max\{0,k- \alpha(d+k) -\zeta(d+k)\}=0.$$ Using the above theorem, we get 
$(\mathbb F_2 \otimes_{\mathcal A} P_5)_{5(2^s-1) + 2^s53} \cong (\mathbb F_2 \otimes_{\mathcal A} P_5)_{5(2^0-1) + 2^053}$ for all $s\geqslant 0.$

Therefore, we need only to study $(\mathbb F_2 \otimes_{\mathcal A} P_5)_{5(2^s-1) + 24.2^{s}}$ for $s=0$ and $s=1.$

\medskip
{\bf Case s=0.}

\medskip
Denote $\omega_{(1)}=(4,4,3), \ \omega_{(2)}=(4,4,1,1), \ \omega_{(3)}=(4,2,2,1), \ \omega_{(4)}=(4,2,4).$ We give a direct summand decomosition of the 
$\mathbb F_2$-vector spaces $(\mathbb F_2 \otimes_{\mathcal A} P_5)_{5(2^0-1) + 24.2^{0}}$ as follows:

\medskip
\begin{thm}\label{mdbd} We have a direct summand decomposition of the $\mathbb F_2$-vector spaces
$$ (\mathbb F_2 \otimes_{\mathcal A} P_5)_{24} =(QP_5^0)_{24}  \oplus QP_5^+(\omega_{(1)}) \oplus QP_5^+(\omega_{(2)}) \oplus  QP_5^+(\omega_{(3)}) \oplus  QP_5^+(\omega_{(4)}).$$
\end{thm}

\medskip
Recall that $(\mathbb{F}_{2}{\otimes}_{\mathcal{A}}P_{4})_{10}$ is an $\mathbb F_2$-vector space of dimension 70 with a basis consisting of all the classes represented by the monomials $w_j, \  1 \leqslant j \leqslant 70$ . Consequently, $|B_4(10)| = 70,$ (see Sum~\cite{su15}).

Since $\mu (24)=4,$ Theorem \ref{dlk} implies that the squaring operation
$$ \widetilde{Sq}_*^0: (\mathbb{F}_{2}{\otimes}_{\mathcal{A}}P_{4})_{24} \longrightarrow (\mathbb{F}_{2}{\otimes}_{\mathcal{A}}P_{4})_{10}$$
is an isomorphism of $GL_4$-module. Hence, $|B_4(24)|=|B_4(10)|=70,$ and therefore the set $[B_4(24)] = \{[\phi_4(w_j)]: w_j \in B_4(10)\}$ is a basis of the $\mathbb F_2$-vector space $(\mathbb{F}_{2}{\otimes}_{\mathcal{A}}P_{4})_{24},$ where $\phi_4(u)=x_1x_2x_3x_4u^2.$

Using Proposition \ref{mdp0}, we obtain $[B_5^0(24)] = \{ [a_t] : \ a_t \in \bigcup _{i=1}^5f_i(B_4(24)), 1\leqslant t \leqslant 350 \}$ is a basis of the $\mathbb F_2$-vector space $(QP_5^0)_{24}.$ Consequently, $\dim (QP_5^0)_{24}=| \bigcup _{i=1}^5f_i(B_4(24))| =350.$

\begin{rem}\label{nx1} We recall a result in Mothebe-Kaelo-Ramatebele \cite{m-k-o16} as follows. 

Set $\mathcal M_{(k,r)}=\{J= (j_1, j_2, \ldots, j_r) : 1 \leqslant j_1 <\ldots < j_r \leqslant k\}$, $1 \leqslant r < k$. For $J\in \mathcal M_{(k,r)}$, define the homomorphism $f_J : P_r \to P_k$ of algebras by substituting $f_J(x_t) = x_{j_t}$ with $1 \leqslant t \leqslant r$. Then, $f_J$ is a monomorphism of $\mathcal A$-modules. We have a direct summand decomposition of the $\mathbb F_2$-vector subspaces: 
$$ QP_k^0  = \bigoplus\limits_{1\leqslant r\leqslant k-1}\bigoplus\limits_{J \in \mathcal M_{(k,r)}}(Qf_J(P_r^+)),$$
where, $Qf_J(P_r^+) = \mathbb{F}_2\otimes_{\mathcal{A}}f_J(P_r^+)$. 

In degree $n,$ we have $\dim (Qf_J(P_r^+))_n = \dim (QP_r^+)_n$ and $|\mathcal M_{(k,r)}| = \binom kr$. Hence, combining with Theorem \ref{dlwo} we get 
\begin{equation*}\label{ctscbs} \dim (QP_k^0)_n = \sum\limits_{\mu(n) \leqslant r\leqslant k-1}\binom{k}{r}\dim(QP_r^+)_n.
\end{equation*} 
And therefore, using the results in Sum \cite{su15}, one gets $\dim (QP_5^0)_{24} =350. $  
\end{rem}

\medskip
\begin{thm}\label{mdbd} 
$$ \dim(QP_5(\omega_{(d)})) = 
\begin{cases}75, &\text{if } \ d=1,\\  
145, & \text{if } \ d=2,\\ 
390, & \text{if } \ d=3,\\ 
1, & \text{if } \ d=4.\\ 
\end{cases}
$$
\end{thm}

The proof of the above theorem is too long and very technical. It is proved by explicitly determining all admissible monomials of $P_5(\omega_{(d)})$ for $d \in \{1, 2, 3, 4\}.$

From the above results, we get the corollary following. 
\medskip
\begin{corl}\label{thm2} There exist exactly $961$ admissible monomials of degree twenty-four in $P_5.$ Consequently, 
$\dim (\mathbb F_2 \otimes_{\mathcal A} P_5)_{24} = 961.$ 
\end{corl}

\medskip
{\bf Case s=1.}

\medskip
Since Kameko's homomorphism $(\widetilde{Sq}^0_*)_{(5,24)}$
is an epimorphism, we have 
$$(\mathbb F_2 \otimes_{\mathcal A} P_5)_{53} \cong \text{Ker}(\widetilde{Sq}^0_*)_{(5,24)} \bigoplus \text{Im}(\widetilde{Sq}^0_*)_{(5,24)}.$$ We have the following

\medskip
\begin{thm}\label{md424} $\text{Im}(\widetilde{Sq}^0_*)_{(5,24)}$ is isomorphic to a subspace of $(\mathbb F_2 \otimes_{\mathcal A} P_5)_{53}$ generated by all the classes represented by the admissible monomials of the form $x_1x_2\ldots x_5u^2,$ for every $u \in B_5(24).$ Consequently, 
$\dim (\text{Im}(\widetilde{Sq}^0_*)_{(5,24)}) = 961.$ 
\end{thm}

Next, we explicitly determine $\text{Ker}(\widetilde{Sq}^0_*)_{(5,24)}$ by giving a direct summand decomposition of the $\mathbb F_2$-vector spaces $\text{Ker}(\widetilde{Sq}^0_*)_{(5,24)}$ as follows.

\begin{thm}\label{md424} We have a direct summand decomposition of the $\mathbb F_2$-vector spaces
$$\text{Ker}(\widetilde{Sq}^0_*)_{(5,24)} =  (QP_5^0)_{53} \oplus QP_5^+(3,3,3,2,1).$$
\end{thm}

\medskip
\begin{rem}\label{nx2} From the result in \cite{m-k-o16} we have

\begin{equation*}\label{ctscbs} \dim (QP_k^0)_{53} = \sum\limits_{\mu(53) \leqslant r\leqslant 4}\binom{5}{r}\dim(QP_r^+)_{53}.
\end{equation*} 

Since $\mu(53)=3,$ $\dim(QP_3^+)_{53} =8$ (see Kameko \cite{ka90}) and $\dim(QP_4^+)_{53} =88$ (see Sum \cite{su15}), one gets $$\dim (QP_5^0)_{53} =\binom{5}{3}.\dim(QP_3^+)_{53} + \binom{5}{4}.\dim(QP_4^+)_{53} =520.$$
And therefore the set $ \{ [b_t] : \ b_t \in \bigcup _{i=1}^5f_i(B_4(53)), 1\leqslant t \leqslant 520 \}$ is a basis of the $\mathbb F_2$-vector space $(QP_5^0)_{53}.$
\end{rem}

We denote by $\mathbb M^d(n) =  \{ [\bigcup _{i=1}^5 x_i^{2^d-1}f_i(x)] : x \in B_{k-1}(n-2^d+1) \}$ and set  
$\mathbb M= \text{Span} \{ u: \ u \in \bigcup _{d=1}^5\mathbb M^d(53) \ \text{and} \ \omega(u)=(3, 3, 3, 2, 1) \}.$ Hence, one gets the theorem following.

\medskip
\begin{thm}\label{md424} The following statements are true:

\medskip
{\rm i)} $\mathbb M$ is the $\mathbb F_2$-vector subspaces of $QP^+_5(\omega),$ where $\omega=(3, 3, 3, 2, 1).$

{\rm ii)} Assume $QP^+_5(\omega)=\mathbb M \oplus \mathbb N.$ Then, $\dim (\mathbb M) = 389$ \  and \ $\dim (\mathbb N).$ 
\end{thm}

The above theorem is proved by explicity determining all admissible monomials in $P^+_5(\omega).$ 

On the other hand, we see that $|B_5(25)|=1240$ (see Sum \cite{su19}) and $|B_5(19)|=905$ (see Tin \cite{tiLA}). Since $\mu(25)=3=\alpha(25+\mu(25)),$ using the result in Sum \cite{su15}, one gets $$|B_6(5(2^{v}-1)+2^{v}25)|=(2^6-1)|B_5(25)|, \   \text{for any integer}\  v \geqslant k-1=5.$$ 
So, we obtain the corollary following.

\medskip
\begin{corl}\label{cor61} There exist exactly $78120$ admissible monomials of degree $m_1=5(2^v -1)+25.2^{v}$ in $P_6,$ for any $v > 4.$ Consequently, $\dim (\mathbb F_2 \otimes_{\mathcal A} P_6)_{m_1} = 78120.$ 
\end{corl}

Similarly, it is easy to check that $\mu(19)=3=\alpha(\mu(19)+19),$ using the result in Sum \cite{su15} we get the following.

\medskip
\begin{corl}\label{cor62} There exist exactly $57015$ admissible monomials of degree $m_2=5(2^r -1)+19.2^{r}$ in $P_6,$ for any $r \geqslant 5.$ Consequently, $\dim (\mathbb F_2 \otimes_{\mathcal A} P_6)_{m_2} = 57015.$ 
\end{corl}

Next, the aim of our paper is to verify this conjecture for the case $k=5$ and the internal degree $5(2^s-1)+24.2^{s},$ with $s$ an arbitrary non-negative integer. We get the following theorem.

\medskip
\begin{thm}\label{thm3} Singer's conjecture is true for $k=5$ and the degree $5(2^s-1)+24.2^{s},$ with $s$ an arbitrary non-negative integer.
\end{thm}

\medskip
\begin{rem}\label{remark} We prove the above theorem by using the admissible monomial basis of degree $5(2^s-1)+24.2^{s}$ in $P_5$ to explicitly compute the vector space $(\mathbb F_2 \otimes_{\mathcal A} P_5)^{GL_5}_{5(2^s-1)+24.2^{s}}$ and combining the computation of  $\text{\rm Ext}_{\mathcal A}^{5, 5+*}(\mathbb F_2, \mathbb F_2)$ by Lin \cite{li08} and Chen \cite{ch93}.
Our method is different from Boardman's method in \cite{bo93} for $k=3$. Recall that Boardman \cite{bo93} computed the space $(\mathbb F_2 \otimes_{\mathcal A} P_3)^{GL_3}$ by using a basis consisting of the all the classes represented by certain polynomials in $P_3.$ It is difficult to use his method for $k \geqslant 4.$ Moreover, our approach can be applied for $k=4$ by using the admissible monomial basis for $P_4$ (see Sum~\cite{su15}). We hope that the conjecture is also true in this case.
\end{rem} 

The proofs of the results of this Note will be pulished in detail elsewhere.

\medskip\noindent {\bf Acknowledgment.} I would like to thank Ho Chi Minh City University of Technology and Education for supporting this work.


{}

\vskip0.7cm

\end{document}